\documentclass[12pt]{article}
\usepackage{amsmath}
\usepackage{graphicx,psfrag,epsf}
\usepackage{enumerate,comment}
\usepackage{natbib}
\usepackage{url} 
\usepackage{hyperref}

\usepackage{latexsym, epsfig, verbatim, color, amssymb, lscape,graphics,natbib}
\usepackage[margin=1in]{geometry}
\newcommand{\blind}{1}

\begin{document}

\date{}

\def\spacingset#1{\renewcommand{\baselinestretch}%
{#1}\small\normalsize} \spacingset{1}


\if1\blind
{
  \title{\bf Fast Algorithm for Calculating Probability of Chess Winning Streaks}
  \author{Guoqing Diao\hspace{.2cm}\\
    Department of Biostatistics and Bioinformatics, George Washington University
    }
  \maketitle
} \fi

\if0\blind
{
  \bigskip
  \bigskip
  \bigskip
  \begin{center}
    {\LARGE\bf Fast Algorithm for Calculating Probability of Chess Winning Streaks}
\end{center}
  \medskip
} \fi

\bigskip
\begin{abstract}Motivated by the controversy in the chess community, where Hikaru Nakamura, a renowned grandmaster, has posted multiple impressive winning streaks over the years on the online platform chess.com, we derive the probabilities of various types of streaks in online chess and/or other sports. Specifically, given the winning/drawing/losing probabilities of individual games, we derive the probabilities of "pure" winning streaks, non-losing streaks, and "in-between" streaks involving at most one draw over the course of games played in a period. The performance of the developed algorithms is examined through numerical studies. 
\end{abstract}

\noindent%
{\it Keywords:} conditional probability, Monte-Carlo simulations, winning streaks.
\vfill

\newpage
\spacingset{1.45} 
\section{Introduction}

Hikaru Nakamura, a famous chess grandmaster currently ranked \#2 in the United States and \#3 in the World, is an active player on chess.com [2], the leading platform for online chess. Over the years, he has played tens of thousands of games on the site, achieving multiple impressive streaks between 2014 and 2024. Notably, Hikaru has posted perfect or near perfect streaks, including 54/54 in 2016, 56.5/57 in 2018, 61/61 in 2020, 52.5/53 in 2020, and 45.5/46 in 2023.  At  first sight, such long winning streaks seem almost impossible. For example, with a 90\% win probability, the chances of winning 61 consecutive games is around 1.6\%. This probability drops to 1.2e-6 if the win rate is 80\%. Indeed, Hikaru’s remarkable streaks have sparked controversy within the chess community. New York Times reported in December 2023 that former chess world champion Vladimir Kramnik "insinuated that Hikaru had probably cheated while playing on the site" [3]. In a blog post on chess.com, Kramnlik wrote:

"Having checked Hikaru's statistics carefully, I have found NUMEROUS low probabilities performances both of him and some of his opponents. Some of which have EXTREMELY low mathematical probability, according to mathematicians. Way below one percent, according to the calculations of those professional mathematicians."

In a recent report, Rosenthal (2024) [4] argued that one should not focus solely on the games within the streaks because these occurred over the course of more than 57,000 online games. Rosenthal (2024) first assessed the probabilities of chess outcomes through chess ratings and the expected scores [4]. A crude approximation (or lower bound) of the probability of winning streaks was then derived by partitioning the total number of $n$ games into $n/k$ non-overlapping sequences, where $k$ represents the number of games in a winning streak. Monte Carlo simulations were also conducted to estimate the probabilities of Hikaru's winning streaks based on 100 simulations. Rosenthal's findings showed that Hikaru's  win streaks are not particularly surprising [4]. A timeline of these events was later reported on chess.com news [5]. 

In this short note, we aim to calculate the exact probabilities of various types of streaks. Given the winning/drawing/losing probabilities of individual games, we derive a closed-form formula for the exact probability that there is at least one winning streak of $k$ or more games out of a total of $n$ games, for any integers $k$ and $n$. We also derive a tight lower bound of the probability of the so-called "in-between" streak involving at most one draw. Numerical studies demonstrate the algorithms based on these derivations perform well.

\section{Algorithms}
Rosenthal (2024) [4] defined three types of streaks: (1) "pure" winning streak; (2) non-losing streak; and (3) "in-between" streak involving at most one draw. In the first case, only wins continue the streak, whereas any draw or loss ends it. In the second case, only a loss ends a streak whereas wins or draws continue it. In the third case, wins or a single draw continue the streak, while a loss or a second draw ends it. 

We start with the first type of streak, i.e., the "pure" winning streak. Suppose there are $n$ games and only two possible outcomes for each game, a win or a loss. The losing probability for the $i$th game is $q_i$, and the winning probability is $1-q_i$. We assume that the outcomes of the $n$ games are mutually independent. We are interested in calculating the probability of no winning streak of $k$ or more consecutive games. We define $p_{m,k}$ as such a probability in the first $m$ games. It is obvious that
\[
p_{m, k} = 1, ~~~m=1,...,k-1,
\]
\[
p_{k,k} = 1 - \prod_{i=1}^k (1-q_i),
\]
and
\[
p_{k, k-1} = 1 - \prod_{i=1}^{k-1} (1-q_i) - \prod_{i=2}^k (1- q_i) - \prod_{i=1}^k (1-q_i).
\]

Suppose we have $p_{1,k},...,p_{m,k}$. We show how to calculate $p_{m+1,k}$. Now consider the $(m+1)$th game. There are two possible outcomes. If the player loses the $(m+1)$th game, then the conditional probability of  no winning streak of $k$ or more games in the first $m+1$ games is $p_{m,k}$. If the player wins the $(m+1)$th game, the conditional probability is $p_{m,k}$ minus the probability of the event that the last $k-1$ games of the first $m$ games are won, and the $(m-(k-1))$th game is lost, and there is no winning streak of $k$ or more games in the first $m-k$ games. We can show that $p_{m+1,k}$ satisfies the following equation:
\[
p_{m+1,k} = q_{m+1} p_{m,k} + (1 - q_{m+1}) \left( p_{m,k} - q_{m-(k-1)} p_{m-k,k}\prod_{j=m-(k-1)+1}^m(1-q_j) \right).
\]

Now, we consider the second type of streak, i.e., the non-losing streak. We can modify the above algorithm for "pure" winning streak by combining the probabilities of winning or drawing a game as the new probability of winning a game and assuming no draws. Then, the algorithm for the "pure" winning streak can be applied directly.

Finally, we consider the third type of streak, i.e., the "in-between" streak involving at most one draw. The probability of the "in-between" streak turns out to be much more complicated than those of the two previous cases. In this case, there are three possible outcomes of each game: win, draw, and loss, resulting in a score of 1, 1/2, and 0, respectively. We denote the probabilities of losing, drawing, and winning the $i$th game by $q_{i,1}, q_{i,2}$, and $q_{i,3}$, respectively. Clearly, $\sum_{j=1}^3 q_{i,j} = 1$ for any $i=1,...,n$. We define $D_{m,k+1/2}$ as the event that there are no  consecutive $k+1$ games out of the first $m$ games such that the player scores $k+1/2$ points or more. Such an "in-between" streak means that the player has $k+1$ wins or $k$ wins and one draw  in at least one of the $k+1$ consecutive games out of the first $m$ games. We define $g_{m, k+1/2}=P(D_{m, k+1/2})$.

Clearly, for any $m \leq k$, $g_{m,k+1/2}=1$. We aim to derive $g_{m, k+1/2}$ for $m\geq k+1$. Now consider $g_{k+1, k+1/2}$. To score $k+1/2$ or more points out of the first $k+1$ games is equivalent to winning all games or winning $k$ games and drawing one game. Hence, 
\[
g_{k+1,k+1/2} = 1 - \prod_{i=1}^{k+1} q_{i,3} - \sum_{i=1}^{k+1} q_{i,2} \left(\prod_{j=1, j\neq i}^{k+1} q_{j,3} \right).
\]
Suppose we have $g_{1, k+1/2},...,g_{m, k+1/2}$. Now consider the $(m+1)$th game. Let $A_i$ be the score of the $i$th game that takes values of 0, 1/2, and 1, corresponding to a loss, a draw, and a win, respectively. If the $(m+1)$th game is a loss, this game does not contribute any additional probability of scoring $k+1/2$ or more points consecutively, and hence $P(D_{m+1, k+1/2} | A_{m+1}=0) = g_{m, k+1/2}$. If the $(m+1)$th game is a draw, the additional probability of scoring $k+1/2$   or more points out of any consecutive $k+1$ games is the probability of the following events: (1) the $k$ games from the $(m-k+1)$th game to the $m$th game are wins; (2) the $(m-k)$th game is a loss; and (3) there is no any sequence of consecutive $k+1$ games out of the first $m-k-1$ games such that the player scores $k+1/2$ points or more (i.e., $D_{m-k-1,k+1/2}$). These events mean that there are no streaks of $k+1/2$ points in the first $m$ games, but there will be a streak in the first $m+1$ games because of the outcome of the $(m+1)$th game.  We have 
\[
P(\bar{D}_{m+1,k+1/2} | A_{m+1}=1/2) = P(\bar{D}_{m,k+1/2}) + 
  P(D_{m-k-1,k+1/2}) q_{m-k,1} \prod_{j=m-k+1}^m q_{j, 3},
\]
where $\bar{D}_{m+1,k+1/2}$ is the complement of $D_{m+1,k+1/2}$. Hence,
\[
P({D}_{m+1,k+1/2} | A_{m+1}=1/2) = g_{m,k+1/2} - 
  g_{m-k-1,k+1/2} q_{m-k,1} \prod_{j=m-k+1}^m q_{j, 3}.
\]

Now consider the last possible scenario for the $(m+1)$th game:  a win. The additional probability of scoring $k+1/2$   or more points out of any consecutive $k+1$ games is the probability of the following events: (1a) the $k$ games from the $(m-k+1)$th game to the $m$th game are wins or $k-1$ games are wins, and the other one is a draw; (2a) the $(m-k)$th game is a loss; and (3a) there is no any sequence of consecutive $k+1$ games out of the first $m-k-1$ games such that the player scores $k+1/2$ points or more (i.e. $D_{m-k-1,k+1/2}$); or (1b) $k-1$ games out of  the $k$ games from the $(m-k+1)$th game to the $m$th game are wins and the other one is a draw; (2b) the $(m-k)$th game is a draw; and (3b) there are no "in-between" streaks with $k+1/2$ points or more out of the first $m$ games. The probability of events (1a)-(3a) is 
\[
  P(D_{m-k-1,k+1/2}) q_{m-k,1} \left\{ \prod_{j=m-k+1}^m q_{j, 3}  + 
  \sum_{i=m-k+1}^m q_{i,2} \left( \prod_{j=m-k+1, j\neq i}^m q_{j, 3}\right) \right\}.
\]
The probability of events (1b)-(3b) takes a more complicated form. We first consider the scenario where the $(m-k+1)$th game is a draw, and the next $k-1$ games are wins. Conditional on this event and (2b), event (3b) is equivalent to $D_{m-k-1, k+1/2}\textbackslash  B_{m-2k-1, m-k-1, k}$, where $B_{a, b, k}$ defines the events such that there are no "in-between" streaks of $k+1/2$ or more points out of the first $b$ games, the $a$th game is a loss (if $b-a=k$) or the $a$th game is not a win (if $b-a<k$), and the next $b-a-1$ games are wins. The corresponding conditional probability of the event (3b) is then 
\[
g_{m-k-1,k+1/2} - h_{m-2k-1,m-k-1,k},
\]
where for $b-a=k$
\[
h_{a,b,k} = 
\begin{cases}
  \prod_{j=a+1}^b q_{j,3} & a=0\\
  q_{a,1} \prod_{j=a+1}^b q_{j,3} & a=1\\
g_{a-1,k+1/2}q_{a,1} \prod_{j=a+1}^b q_{j,3} & a>1\\
\end{cases},
\]
and for $b-a<k$
\[
h_{a,b,k} 
\begin{cases}
  = \prod_{j=a+1}^b q_{j,3} & a=0\\
  = (q_{a,1}+q_{a,2}) \prod_{j=a+1}^b q_{j,3} & a=1\\
 \approx g_{a-1,k+1/2}(q_{a,1}+q_{a,2}) \prod_{j=a+1}^b q_{j,3} & a>1\\
\end{cases}.
\]
Note that when $b-a<k$ and $a>1$, the probability of $B_{a,b,k}$ is complicated, and we resort to an approximation, which is also an upper bound of $h_{a,b,k}$. This approximation, as shown in the numerical studies, has minimal impact on the accuracy  of $g_{m+1,k+1/2}$. 

Now consider the general scenario such that the $l$th ($m-k+1 \leq l \leq m$) game is a draw, and the other $k-1$ games from the $(m-k+1)$th game to the $m$th game are wins. Then conditional on this event and (2b), event (3b) is equivalent to $D_{m-k-1, k+1/2}\textbackslash \left(\cup_{a=m-2k-1}^{l-k-2} B_{a, m-k-1} \right)$. The corresponding conditional probability of (3b) is 
\[
g_{m-k-1,k+1/2} - \sum_{a=m-2k-1}^{l-k-2} h_{a, m-k-1,k}. 
\]

The above derivations  suggest that
\[
\begin{split}
P(\bar{D}_{m+1,k+1/2} | A_{m+1}=1) = & P(\bar{D}_{m,k+1/2}) + 
  P(D_{m-k-1,k+1/2}) q_{m-k,1} \\
  & \times \left\{ \prod_{j=m-k+1}^m q_{j, 3}  + 
  \sum_{i=m-k+1}^m q_{i,2} \left( \prod_{j=m-k+1, j\neq i}^m q_{j, 3}\right) \right\} \\
  & + \sum_{l=m-k+1}^m \left(g_{m-k-1,k+1/2} - \sum_{a=m-2k-1}^{l-k-2} h_{a, m-k-1,k} \right) \\
  & \times q_{m-k,2} q_{l,2} \prod_{j=m-k+1, j\neq l}^m q_{j,3}.\\
\end{split}
\] Hence, we have
\[
\begin{split}
P({D}_{m+1,k+1/2} | A_{m+1}=1) = & g_{m,k+1/2} - 
  g_{m-k-1,k+1/2} q_{m-k,1} \\
  & \times \left\{ \prod_{j=m-k+1}^m q_{j, 3}  + 
  \sum_{i=m-k+1}^m q_{i,2} \left( \prod_{j=m-k+1, j\neq i}^m q_{j, 3}\right) \right\}\\
  & - \sum_{l=m-k+1}^m \left(g_{m-k-1,k+1/2} - \sum_{a=m-2k-1}^{l-k-2} h_{a, m-k-1,k} \right) \\
  & \times q_{m-k,2} q_{l,2} \prod_{j=m-k+1, j\neq l}^m q_{j,3}.\\
\end{split}
\]
Finally, we obtain
\[
g_{m+1,k+1/2} = q_{m+1,1} g_{m,k+1/2} + q_{m+1,2} P({D}_{m+1,k+1/2} | A_{m+1}=1/2) + q_{m+1,3} P({D}_{m+1,k+1/2} | A_{m+1}=1).
\]

\section{Numerical Studies}
We conduct numerical studies to examine the performance of the proposed algorithms. In the first set of numerical studies, we consider 30,000 games and four scenarios for the distribution of the winning probabilities of individual games: (1) $U(0.85, 1)$; (2) $U(0.8, 1)$; (3) $U(0.75, 1)$; and (4) $U(0.7, 1)$. Figure 1 presents the probabilities of having "pure" winning streaks of $k$ or more games $(k=1,...,150)$ out of 30,000 games. As expected, the higher the winning probabilities of individual games, the higher the chance of having long winning streaks. For example, when the winning probabilities of individual games follow $U(0.85, 1)$, the probabilities of having winning streaks of 50+, 75+, 100+, 125+, and 150+ games are 1, 0.999, 0.605, 0.124, and 0.019, respectively. To check the accuracy of the algorithm, we also conduct simulations to estimate the probabilities of winning streaks by generating 10,000 random sequences of 30,000 games. Table 1 presents the true and estimated probabilities of having winning streaks, where $k$ is chosen such that the probabilities of having winning streaks of $k$ or more games are close to 0.90, 0.75, 0.50, 0.25, and 0.10. In all cases, the estimated probabilities of "pure" winning streaks are close to the true probabilities calculated by the proposed algorithm.

In the second set of numerical studies, we consider 10,000 games and the "in-between" streaks. We generate the drawing probabilities from $U(0,0.2)$, whereas four scenarios for the distribution of the losing probabilities are considered: (1) $U(0, 0.15)$; (2) $U(0, 0.2)$; (3) $U(0, 0.25)$; and (4) $U(0, 0.3)$.  Figure 2 presents the probabilities of having "in-between" streaks of $k+1/2$ or more  $(k=1,...,60)$ out of 10,000 games. Again, as expected, the lower the losing probabilities of individual games, the higher the chance of having long "in-between" streaks. For example, when the losing probabilities of individual games follow $U(0, 0.15)$, the probabilities of having "in-between" streaks of (15+1/2)+, (25+1/2)+, (35+1/2)+, (45+1/2)+, and (55+1/2)+ games are 1, 1, 1, 0.763, and 0.225, respectively. As in the first set of studies, we also conduct simulations to estimate the probabilities of  "in-between" streaks allowing for at most one draw by generating 10,000 random sequences of 10,000 games. Table 2 presents the true and estimated probabilities of having such "in-between" streaks. Although in the algorithm, we use an approximation to $h_{a,b,k}$ when $b-a<k$ and $a>1$, the estimated probabilities are still very close to the probabilities calculated from the proposed algorithm, suggesting that the approximation is reasonably accurate. 

We implement the proposed algorithms in a C program. The C program is computationally efficient, and it takes less than one second on a Mac Studio with the Apple M2 Ultra Chip to calculate $n\times k$ probabilities of winning streaks in each setting of the above numerical studies. Note that for the first set of studies, $n=30,000$ and $k=150$, and for the second set of studies, $n=10,000$ and $k=60$.

\begin{figure}[t]
        \label{fig1} 
         \caption{Probabilities of having "pure" winning streaks of $k$ or more games out of a total of 30,000 games with winning probabilities of individual games from: (1) $U(0.85, 1)$; (2) $U(0.8, 1)$; (3) $U(0.75, 1)$; and (4) $U(0.7, 1)$. }
        \centering\includegraphics[width=1.0\textwidth]{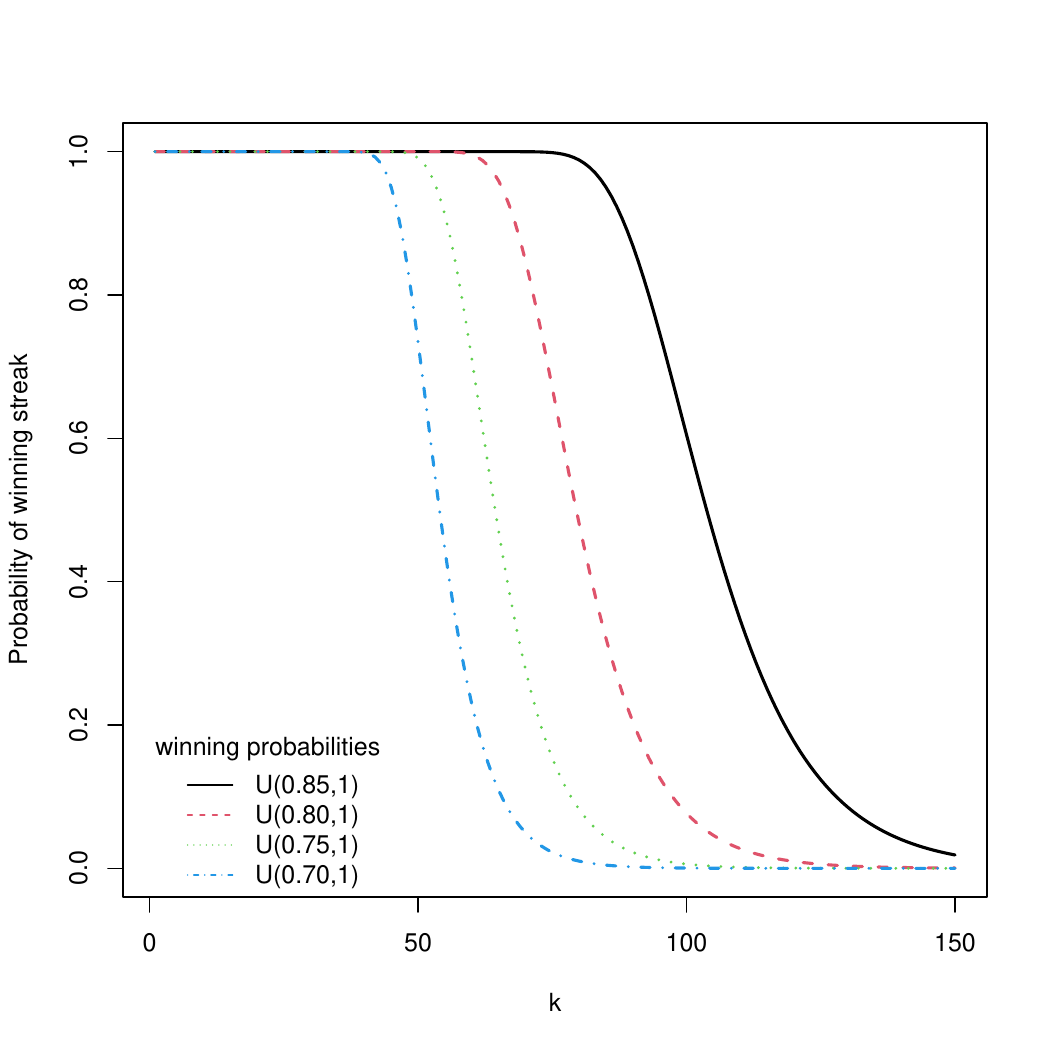}
        
\end{figure}

\begin{figure}[t]
        \label{fig2} 
         \caption{Probabilities of having "in-between" streaks of $k+1/2$ (i.e., allowing for at most one draw) or more games out of a total of 10,000 games with drawing probabilities of individual games from $U(0, 0.2)$ and losing probabilities following: (1) $U(0, 0.15)$; (2) $U(0, 0.2)$; (3) $U(0, 0.25)$; and (4) $U(0, 0.3)$. }
        \centering\includegraphics[width=1.0\textwidth]{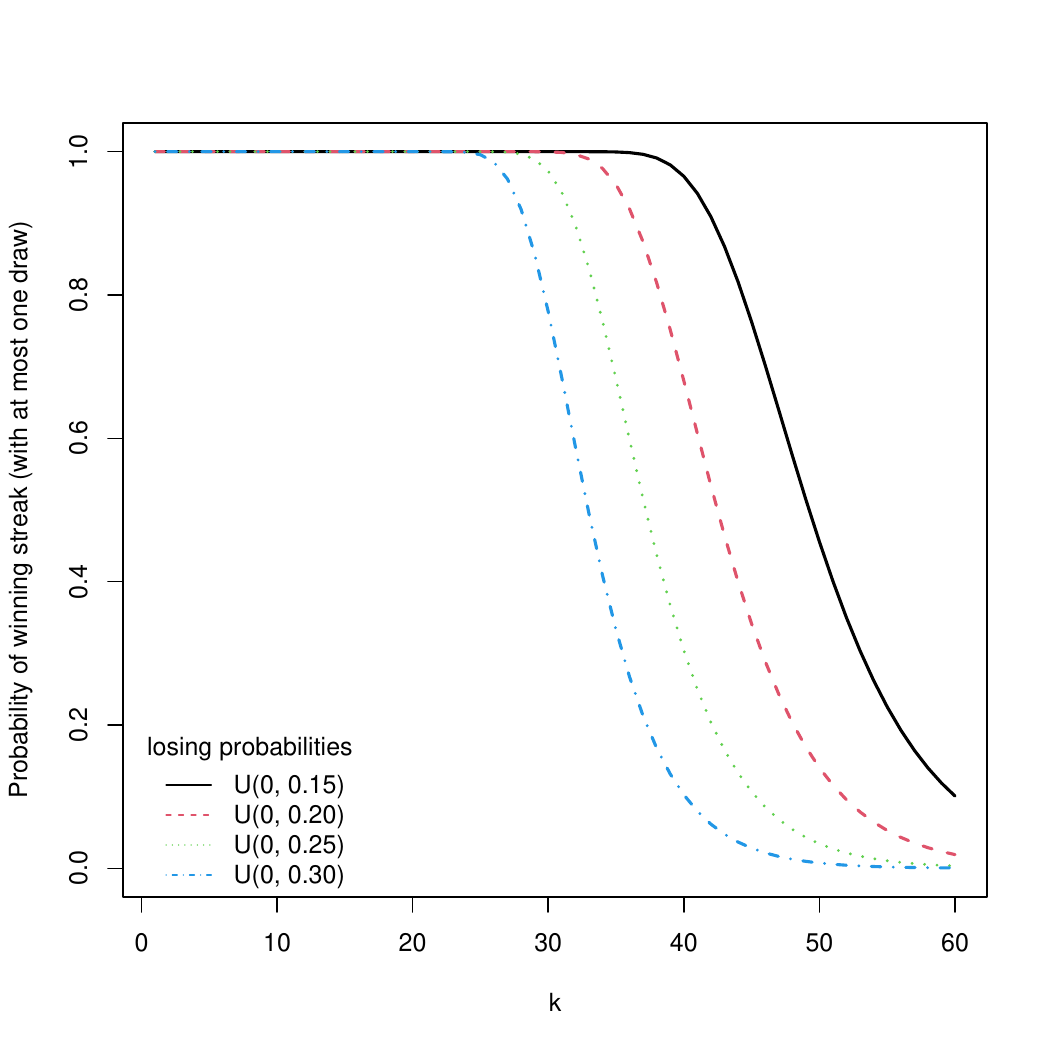}
        
\end{figure}

\begin{table}

\caption{True and estimated probabilities (based on 10,000 simulations) of having "pure" winning streaks of $k$ or more games out of 30,000 games with winning probabilities of individual games from: (1) $U(0.85, 1)$; (2) $U(0.8, 1)$; (3) $U(0.75, 1)$; and (4) $U(0.7, 1)$}
\begin{center}
\begin{tabular}{lcccclccc}\hline
$k$ &True & Est. & &$k$ &True & Est.\\
\multicolumn{3}{c}{$U(0.85,1)$} & &\multicolumn{3}{c}{$U(0.80,1)$}\\
88 & 0.9069 & 0.9067 && 68 & 0.9024 & 0.9003\\
94 & 0.7735 & 0.7765 && 72 & 0.7820 & 0.7828\\
103 & 0.5203 & 0.5232 && 79 & 0.5170 & 0.5174\\
115 & 0.2502 & 0.2616 && 87 & 0.2690 & 0.2740\\
127 & 0.1068 & 0.1060 && 97 & 0.1037 & 0.0991\\
\multicolumn{3}{c}{$U(0.75,1)$} & &\multicolumn{3}{c}{$U(0.70,1)$}\\
55 & 0.9117 & 0.9164 && 46 & 0.9217 & 0.9213\\
59 & 0.7582 & 0.7559 && 49 & 0.7896 & 0.7886\\
64 & 0.5164 & 0.5178 && 53 & 0.5566 & 0.5636\\
70 & 0.2777 & 0.2771 && 59 & 0.2639 & 0.2599\\
78 & 0.1059 & 0.1058 && 65 & 0.1090 & 0.1053\\
\hline
\end{tabular}
\end{center}
\end{table}

\begin{table}

\caption{True and estimated probabilities (based on 10,000 simulations) of having "in-between" streaks of $k+1/2$ (i.e., allowing for at most one draw) or more games out of 10,000 games with drawing probabilities of individual games from $U(0, 0.2)$ and losing probabilities following: (1) $U(0, 0.15)$; (2) $U(0, 0.2)$; (3) $U(0, 0.25)$; and (4) $U(0, 0.3)$}
\begin{center}
\begin{tabular}{lcccclccc}\hline
$k$ &True & Est. & &$k$ &True & Est.\\
\multicolumn{3}{c}{$U(0, 0.15)$} & &\multicolumn{3}{c}{$U(0, 0.20)$}\\
42 & 0.9091 & 0.9112 && 36 & 0.9200 & 0.9203\\
45 & 0.7625 & 0.7658 && 39 & 0.7515 & 0.7605\\
49 & 0.5146 & 0.5136 && 42 & 0.5338 & 0.5365\\
54 & 0.2621 & 0.2633 && 46 & 0.2889 & 0.2971\\
60 & 0.1008 & 0.1033 && 51 & 0.1160 & 0.1170\\
\multicolumn{3}{c}{$U(0, 0.25)$} & &\multicolumn{3}{c}{$U(0, 0.30)$}\\
31 & 0.9433 & 0.9412 && 28 & 0.9195 & 0.9169\\
34 & 0.7641 & 0.7680 && 30 & 0.7768 & 0.7751\\
37 & 0.5152 & 0.5125 && 32 & 0.5886 & 0.5926\\
40 & 0.3041 & 0.3039 && 36 & 0.2670 & 0.2711\\
45 & 0.1069 & 0.1117 && 40 & 0.1026 & 0.1028\\
\hline
\end{tabular}
\end{center}
\end{table}

\section{Limitations}
We have developed efficient algorithms to calculate the probability of various types of streaks, including the "pure" winning streak, the non-losing streak, and the "in-between" streak involving at most one draw. Numerical studies demonstrated that these algorithms perform very well. 
There are, however, several limitations. First, we use an approximation to $h_{a,b,k}$ for $b-a<k$ and $a>1$. We have to resort to the approximation because it is challenging to exhaust all possible scenarios for $B_{a,b,k}$. On the other hand, this approximation provides an upper bound to $h_{a,b,k}$, and consequently, the calculated probability of having "in-between" streaks is a lower bound. In addition, the impact of this approximation is minimal, as demonstrated by the numerical studies.

A second limitation of the method is that one needs to assume independence between games. This independence assumption is violated in the presence of "hot hands." However, as pointed out by Rosenthal (2024), any persistence or "hot hand" assumption would instead make the streaks more likely [4]. This implies that the probability of streaks calculated under the independence assumption again provides a lower bound for the probability of streaks in the presence of "hot hands."

Another limitation is that we calculate the probabilities of streaks over many games, regardless of the levels of opponents (or the winning/drawing/losing probabilities of individual games). Apparently, having a streak against strong opponents is more challenging than having a streak against weaker opponents. One possible solution to this limitation is to focus on games against players with ratings above a certain threshold.  

\section{Applications to other fields}
Although this work was motivated by Hikaru's impressive winning streaks in online chess, the algorithm can be applied to other online or offline games, e.g., poker, video games, etc. Given reasonable estimates of the games' winning/drawing/losing probabilities, the proposed  algorithms provide an efficient way to calculate the probability of various types of streaks of a player or team in their respective sport. 

\section*{Disclaimer}
This work was not sponsored by chess.com or any other entity. It was done out of the author's free time. 

\section*{Reference:}

[1] International Chess Federation (2025)
\url{https://ratings.fide.com/profile/2016192}, accessed on February 25, 2025.

\noindent [2] \url{https://www.chess.com/players/hikaru-nakamura}

\noindent [3] \url{https://www.nytimes.com/2023/12/25/crosswords/chess-hikaru-vladmir-kramnik-cheating.html}

\noindent [4] Jeﬀrey S. Rosenthal (2024). 
\url{https://probability.ca/jeff/ftpdir/chessstreakpaper.pdf}, accessed on February 25, 2025.

\noindent [5] Tarjei J. Svensen (2024). \url{https://www.chess.com/news/view/nakamura-winning-streaks-statistically-normal-professor-says}

\end{document}